\documentclass[reqno]{amsart}
\usepackage{amssymb,amsmath,amsfonts,amscd,amsthm,pb-diagram}
\usepackage{amsbsy,bm}

\newtheorem{theorem}{Theorem}[section]

\newtheorem{definition}{Definiton}[section]

\newtheorem{lemma}{Lemma}[section]

\theoremstyle{remark}
\newtheorem{remark}{Remark}[section]
\newtheorem{example}{Example}[section]

\begin{document}
\title[Some rigidity phenomenons for Lagrangian submanifolds ]
{On energy gap phenomena of the Whitney spheres in $\mathbb{C}^n$ or $\mathbb{CP}^n$}

\author{Yong Luo, Liuyang Zhang}
\address{Mathematical Science Research Center, Chongqing University of Technology, Chongqing, 400054, China}
\email{yongluo-math@cqut.edu.cn}
\address{Mathematical Science Research Center, Chongqing University of Technology, Chongqing, 400054, China}
\email{13320241808@163.com}

\thanks{The first author was supported by NSF of China, Grant Number 11501421.}

\date{}


\begin{abstract}
 In  \cite{Zh} \cite{LY} Zhang, Luo and Yin initiated the study of Lagrangian submanifolds satisfying ${\rm \nabla^*} T=0$ or ${\rm \nabla^*\nabla^*}T=0$  in  $\mathbb{C}^n$ or $\mathbb{CP}^n$, where $T ={\rm \nabla^*}\tilde{h}$ and $\tilde{h}$ is the Lagrangian trace-free second fundamental form. They proved several rigidity theorems for Lagrangian surfaces satisfying ${\rm \nabla^*} T=0$ or ${\rm \nabla^*\nabla^*}T=0$ in $\mathbb{C}^2$ under proper small energy assumption and gave new characterization of the Whitney spheres in $\mathbb{C}^2$. In this paper we extend these results to Lagrangian submanifolds in $\mathbb{C}^n$ of dimension $n\geq3$  and to Lagrangian submanifolds in $\mathbb{CP}^n$.
\end{abstract}

\maketitle

\numberwithin{equation}{section}
\section{Introduction}\label{sect:1}
Assume that $N^n(4c)$ is the the standard complex space form with standard complex structure $J$, K\"ahler form $\omega$ and metric $\langle , \rangle$, i.e. $N^n(0)=\mathbb{C}^n$ and $N^n(4)=\mathbb{CP}^n$. A real $n$-dimensional submanifold of $N^n(4c)$ is a Lagrangian submanifold if $J$ is an isometric map between its tangent bundle and normal bundle. The most canonical and important examples of Lagrangian submanifolds of $\mathbb{C}^n$ or $\mathbb{CP}^n$ are the Lagrangian subspaces and  Whitney spheres. The Whitney spheres in $\mathbb{C}^n$ are defined by (cf. \cite{We})
\begin{example}
\begin{eqnarray*}
\phi_{r,A}: {\mathbb{S}^n}&{\to}& {\mathbb{C}^n}
 \\(x_1,...,x_{n+1})&\mapsto& \frac{r}{1+x_{n+1}^2}(x_1,x_1x_{n+1},...,x_n,x_nx_{n+1})+A,
\end{eqnarray*}
where $\mathbb{S}^n=\{(x_1,...,x_{n+1})\in {\mathbb{R}^{n+1}}|x_1^2+...+x_{n+1}^2=1\},$ $r$ is a positive number and $A$ is a vector of $\mathbb{C}^n$.
\end{example}
The Whitney spheres in $\mathbb{CP}^n$ are defined by (cf. \cite{CU1}\cite{Ch})
\begin{example}
\begin{eqnarray*}
\phi_\theta: {\mathbb{S}^n}&{\to}& {\mathbb{CP}^n}, \theta>0
 \\(x_1,...,x_{n+1})&\mapsto& [(\frac{(x_1,...,x_n)}{ch\theta+ish\theta x_{n+1}}, \frac{sh\theta ch\theta(1+x^2_{n+1})+ix_{n+1}}{ch^2_\theta+sh^2_\theta x^2_{n+1}})],
\end{eqnarray*}
where $\mathbb{S}^n=\{(x_1,...,x_{n+1})\in {\mathbb{R}^{n+1}}|x_1^2+...+x_{n+1}^2=1\}.$
\end{example}
The Lagrangian subspaces and Whitney spheres $\phi_{r,A}$ in $\mathbb{C}^n$ or the real projective space $\mathbb{RP}^n$ and Whitney spheres $\phi_\theta$ in $\mathbb{CP}^n$ play a similar role with that of totally umbilical hypersurfaces in a real Euclidean space $\mathbb{R}^{n+1}$ or in the unit sphere $\mathbb{S}^{n+1}$, and they are locally characterized by vanishing of the following so called Lagrangian trace free second fundamental form (cf. \cite{CU}\cite{RU}\cite{CU1}\cite{Ch})
\begin{eqnarray}
\tilde{h}(V,W):=h(V,W)-\frac{n}{n+2}\{\langle V, W\rangle H+\langle JV, H\rangle JW+\langle JW, H\rangle JV\},
\end{eqnarray}
where $h$ denotes the second fundamental form and $H=\frac{1}{n}h$ denotes the mean curvature vector field respectively.

Various characterizations of the Lagrangian subspaces, $\mathbb{RP}^n$ or Whitney spheres in $\mathbb{C}^n$ or $\mathbb{CP}^n$ were obtained in \cite{BCM}\cite{CU, CU2,CMU}\cite{CD}\cite{Ch1,Ch2}\cite{LV}\cite{LW}\cite{RU}. In particular, Castro, Montealegre, Ros and Urbano \cite{CU}\cite{RU}\cite{CMU} introduced and studied Lagrangian submanifolds with conformal Maslov form in $\mathbb{C}^n$ or $\mathbb{CP}^n$,  that is Lagrangian submanifolds in $\mathbb{C}^n$ or $\mathbb{CP}^n$ with the 2-form $T=0$, where in local orthonormal basis
\begin{eqnarray}
T_{ij}:=\tfrac1n\sum_m\tilde{h}^{m^*}_{ij,m}=\tfrac{1}{n+2}\big(nH^{i^*}_{,j}-\sum_mH^{m^*}_{,m}\,g_{ij}\big).
\end{eqnarray}
They proved that the only compact(nonminimal) Lagrangian submanifolds in $\mathbb{C}^n$ or $\mathbb{CP}^n$ with  conformal Maslov form (i.e. $T=0$) and null first Betti number are the Whitney spheres. The Whitney spheres in $\mathbb{C}^n$ also play an important role in the study of Lagrangian mean curvature flow \cite{CLM}\cite{SS}.

Recently, Zhang, Luo and Yin \cite{Zh}\cite{LY} initiated the study of Lagrangian submanifolds in $\mathbb{C}^n$ or $\mathbb{CP}^n$ satisfying $\nabla^*T=0$ or $\nabla^*\nabla^*T=0$. In particular, they proved the following results:
\begin{theorem}[\cite{Zh}]
 Assume that $\Sigma\hookrightarrow\mathbb{C}^2$ is a properly immersed complete Lagrangian surface satisfying $\nabla^*T=0$. Then there exists a constant $\epsilon_0>0$ such that if
 $$\int_\Sigma|\tilde{h}|^2d\mu\leq\epsilon_0\ and\ \lim_{R\to+\infty}\frac{1}{R^2}\int_{\Sigma_R}|h|^2d\mu=0,$$
 where $\Sigma_R:=\Sigma\cap B_R(0)$ and $B_R(0)$ denotes the ball centered at $0$ in $\mathbb{C}^2$ with radius $R$, then $\Sigma$ is either a Lagrangian plane or a 2-dimensional Whitney sphere.
\end{theorem}
\begin{remark}
Though it was assumed properness in the above theorem, we see  from the proof in \cite{Zh} that we only need assume that $\Sigma$ is complete.
\end{remark}
\begin{theorem}[\cite{LY}]
 Assume that $\Sigma\hookrightarrow\mathbb{C}^2$ is a Lagrangian sphere satisfying $\nabla^*\nabla^*T=0$. Then there exists a constant $\epsilon_0>0$ such that if
 $$\int_\Sigma|\tilde{h}|^2d\mu\leq\epsilon_0,$$
 then $\Sigma$ is a 2-dimensional Whitney sphere.
\end{theorem}
The aim of this paper is to extend the above results to higher dimensional Lagrangian submanifolds in $\mathbb{C}^n$ and to Lagrangian submanifolds in $\mathbb{CP}^n$. In fact we have
\begin{theorem}\label{thm:1.3}
 Assume that $M^n\hookrightarrow\mathbb{C}^n(n\geq3)$ is a complete Lagrangian submanifold. We have
 \\(i) if $M^n$ satisfies $\nabla^*T=0$, then there exists a constant $\epsilon_0>0$ such that if
 $$\int_M|\tilde{h}|^nd\mu\leq\epsilon_0\ and\ \lim_{R\to+\infty}\frac{1}{R^2}\int_{M_R}|h|^2d\mu=0,$$
 where $M_R$ denotes the geodesic ball in $M^n$ with radius $R$, then  $M^n$ is either a Lagrangian subspace or a Whitney sphere;
 \\
 \\(ii) if $M^n$ is a Lagrangian sphere satisfying $\nabla^*\nabla^*T=0$, then there exists a constant $\epsilon_0>0$ such that if
 $$\int_M|\tilde{h}|^nd\mu\leq\epsilon_0,$$
 then $M^n$ is a Whitney sphere.
\end{theorem}

We would like to point out that compared with the 2-dimensional case, the proof of Theorem \ref{thm:1.3} is much more complicated. Firstly, in the 2-dimensional case we just need to test over a simple Simons' type identity, but in the case of dimension $n\geq3$ we need to estimates the nonlinear terms in a much more complicated Simons' type equality to get a Simons' type inequality (cf. (\ref{eqn:3.17})) and then test over it. Secondly, in the higher dimensional case we need to adapt the original Michael-Simon inequality a little bit to get (\ref{eqn:3.20}) and use it to absorb the "bad term" at the right hand of (\ref{eqn:3.17}).

Similarly, for Lagrangian submanfiods in $\mathbb{CP}^n$ we have
\begin{theorem}\label{thm:1.4}
 Assume that $M^n\hookrightarrow\mathbb{CP}^n(n\geq2)$ is a complete Lagrangian submanifold. We have
 \\(i) if $M^n$ satisfies $\nabla^*T=0$, then there exists a constant $\epsilon_0>0$ such that if
 $$\int_M|\tilde{h}|^nd\mu\leq\epsilon_0 \ and\ \lim_{R\to+\infty}\frac{1}{R^2}\int_{M_R}|h|^2d\mu=0,$$
 where $M_R$ denotes the geodesic ball in $M^n$ with radius $R$,
 then  $M^n$ is the real projective space $\mathbb{RP}^n$ or  a Whitney sphere;
 \\
 \\(ii) if $M^n$ is a Lagrangian sphere satisfying $\nabla^*\nabla^*T=0$, then there exists a constant $\epsilon_0>0$ such that if
 $$\int_M|\tilde{h}|^nd\mu\leq\epsilon_0,$$
 then $M^n$ is the real projective space $\mathbb{RP}^n$ or  a Whitney sphere.
\end{theorem}

Note that similar $L^\frac n2$ pinching theorems for minimal submanifolds in a unit sphere were initiated by Shen \cite{Sh}, and later investigated by Wang \cite{Wa}, Lin and Xia \cite{LX}. $L^\frac n2$ pinching theorems for minimal submanifolds in a Euclidean space was investigated by Ni \cite{Ni}. Generalizations of  $L^\frac n2$ pinching theorems  to submanifolds with parallel mean curvature vector field in a sphere or in a Euclidean space were obtained by Xu \cite{Xu} and Xu and Gu \cite{XG}. Our results could be seem as extensions of their results to more general submanifolds in the Lagrangian setting.

The rest of this paper is organized as follows. In section 2 we give some preliminaries on Lagrangian submanifolds in $N^n(4c)$. In section 3 we prove a Simons' type inequality for Lagrangian submanifolds in $N^n(4c)$, which plays a crucial role in the proof of Theorems \ref{thm:1.3}, \ref{thm:1.4}. Theorem \ref{thm:1.3} is proved in section 4 and Theorem \ref{thm:1.4} is proved in section 5.
\section{Preliminaries}\label{sect:2}~
In this section we collect some basic formulas and results of the Lagrangian submanifolds in a complex space form (cf. \cite{Bl}\cite{Ca}).

Let $N^n(4c)$ be a complete, simply connected, $n$-dimensional K\"ahler manifold
with constant holomorphic sectional curvature $4c$. Let $M^n$ be an $n$-dimensional
Lagrangian submanifolds in $N^n(4c)$. We  denote also by $g$ the metric on
$M^n$. Let $\nabla$ (resp. $\bar\nabla$) be the Levi-Civita
connection of $M^n$ (resp. $N^{n}(4c)$). The Gauss and Weingarten
formulas of $M^n\hookrightarrow N^{n}(4c)$ are given,
respectively, by
\begin{equation}\label{eqn:2.1}
\bar\nabla_XY=\nabla_XY+h(X,Y)\ \ {\rm and}\ \
\bar\nabla_XV=-A_VX+\nabla^{\bot}_XV,
\end{equation}
where $X,Y\in TM^n$ are tangent vector fields, $V\in T^\perp M^n$ is a
normal vector field; $\nabla^{\bot}$ is the normal connection in the
normal bundle $T^\perp M^n$; $h$ is the second fundamental form and
$A_V$ is the shape operator with respect to $V$. From
\eqref{eqn:2.1}, we easily get
\begin{equation}\label{eqn:2.2}
g(h(X,Y),V)=g(A_VX,Y).
\end{equation}
The mean curvature vector $H$ of $M^n$ is defined by $H=\tfrac1n{\rm trace}\,h$.

For Lagrangian submanifolds in $N^n(4c)$, we have
\begin{equation}\label{eqn:2.3}
 \nabla_XJY=J\nabla^\bot_XY,
\end{equation}
\begin{equation}\label{eqn:2.4}
A_{JX}Y=-Jh(X,Y)=A_{JY}X.
\end{equation}
The above formulas immediately imply that $g(h(X,Y),JZ)$ is totally symmetric.

To utilize the moving frame method, we will use the following range
convention of indices:
\begin{equation*}
\begin{gathered}
i,j,k,l,m,p,s=1,\ldots,n; \ \
i^*=i+n \ etc..
\end{gathered}
\end{equation*}

Now, we choose a local {\it adapted Lagrangian frame} $\{e_1,\ldots, e_n,
e_{1^*},\ldots, e_{n^*}\}$ in $N^{n}(4c)$ in such a
way that, restricted  to $M^n$, $\{e_1,\ldots, e_n\}$ is an
orthonormal frame of $M^n$, and $\{e_{1^*}= Je_1,\ldots,e_{n^*}=J e_n\}$ is a orthonormal
frame of $M^n\hookrightarrow N^{n}(4c)$. Let
$\{\theta_1,\ldots,\theta_n\}$ be the dual frame of $\{e_1,\ldots ,
e_n\}$. Let $\theta_{ij}$ and $\theta_{i^*j^*}$ denote the
connection $1$-forms of $TM^n$ and $T^\perp M^n$, respectively.

Put $h_{ij}^{k^*} = g(h(e_i, e_j),J e_k)$. It is easily seen
that
\begin{equation}\label{eqn:2.5}
h_{ij}^{k^*}=h_{ik}^{j^*}=h_{jk}^{i^*},\ \ \forall\ i,j,k.
\end{equation}

Denote by $R_{ijkl}:=g\big(R(e_i,e_j)e_l,e_k\big)$ and
$R_{ijk^*l^*}:=g\big(R(e_i,e_j)e_{l^*},e_{k^*}\big)$
the components of the curvature tensors of $\nabla$ and
$\nabla^{\bot}$ with respect to the adapted Lagrangian frame, respectively.
Then, we get the Gauss, Ricci and Codazzi equations:
\begin{equation}\label{eqn:2.6}
R_{ijkl}=c(\delta_{ik}\delta_{jl}-\delta_{il}\delta_{jk})
+\sum_{m}(h_{ik}^{m^*}h_{jl}^{m^*}-h_{il}^{m^*}h_{jk}^{m^*}),
\end{equation}
\begin{equation}\label{eqn:2.7}
R_{ijk^*l^*}=c(\delta_{ik}\delta_{jl}-\delta_{il}\delta_{jk})
+\sum_{m}(h_{ik}^{m^*}h_{jl}^{m^*}-h_{il}^{m^*}h_{jk}^{m^*}),
\end{equation}
\begin{equation}\label{eqn:2.8}
h^{m^*}_{ij,k}=h^{m^*}_{ik,j},
\end{equation}
where $h^{m^*}_{ij,k}$ is the components of the covariant
differentiation of $h$, defined by
\begin{equation}\label{eqn:2.9}
\sum_{l=1}^nh^{m^*}_{ij,l}\theta_l:=dh_{ij}^{m^*}+\sum_{l=1}^nh^{m^*}_{il}\theta_{lj}
+\sum_{l=1}^nh^{m^*}_{jl}\theta_{li}
+\sum_{l=1}^{n}h^{l^*}_{ij}\theta_{l^*m^*},
\end{equation}

Then from \eqref{eqn:2.5} and \eqref{eqn:2.8}, we have
\begin{equation}\label{eqn:2.10}
h^{m^*}_{ij,k}=h^{i^*}_{jk,m}=h^{j^*}_{km,i}=h^{k^*}_{mi,j}.
\end{equation}

We also have Ricci identity
\begin{equation}\label{eqn:2.11}
 h^{m^*}_{ij,lp}-h^{m^*}_{ij,pl}
=\sum_{k=1}^nh^{m^*}_{kj}R_{kilp}+\sum_{k=1}^nh^{m^*}_{ik}R_{kjlp}
+\sum_{k=1}^nh^{k^*}_{ij}R_{k^*m^{*}lp},
\end{equation}
where $h^{m^*}_{ij,lp}$ is defined by
\begin{equation*}
\sum_p
h^{m^*}_{ij,lp}\theta_p=dh^{m^*}_{ij,l}+\sum_ph^{m^*}_{pj,l}\theta_{pi}
+\sum_ph^{m^*}_{ip,l}\theta_{pj}+
\sum_ph^{m^*}_{ij,p}\theta_{pl}+\sum_{p}h^{p^*}_{ij,l}\theta_{p^*m^*}.
\end{equation*}
The mean curvature vector $H$ of
$M^n\hookrightarrow N^{n}(4c)$ is
$$
H=\tfrac1n\sum\limits_{i=1}^nh(e_i,e_i)=\sum\limits_{k=1}^nH^{k^*}e_{k^*},\
\ H^{k^*}=\tfrac1n\sum_ih^{k^*}_{ii}.
$$
Letting $i = j$ in \eqref{eqn:2.9} and carrying out summation over $i$, we have
\begin{equation*}
H^{k^*}_{,l}\theta_l=dH^{k^*}+\sum_{l}H^{l^*}\theta_{l^*k^*},
\end{equation*}
and we further have
\begin{eqnarray}\label{eqn:12}
{H}^{k^*}_{,i}=H^{i^*}_{,k}
\end{eqnarray} for any $i,k$.

\section{A Simons' type inequality}\label{sect:3}

In this section, inspired by Chao and Dong \cite{CD}, we will derive a new Simons' type inequality for Lagrangian submanifolds in $ N^{n}(4c)$. 

We assume that $M^n \hookrightarrow N^{n}(4c)$ is a Lagrangian submanifold and $n\geq2$, where $N^n(4c)$ is the the standard complex space form of constant holomorphic sectional curvature $4c$ with standard complex structure $J$, K\"ahler form $\omega$ and metric $\langle , \rangle$.

Firstly, we define a trace-free tensor $\tilde h(X,Y)$ defined by
\begin{equation}\label{eqn:3.1}
\tilde h(X,Y)=h(X,Y)-\frac{n}{n+2}\big\{g(X,Y) H+g(J X,H)J Y+g(J Y,H)J X \big\}
\end{equation}
for any tangent vector fields $X,Y$ on $M^n$.

With respect to Lagrangian frame
$\{e_1,\ldots, e_n,
e_{1^*},\ldots, e_{n^*}\}$ in $ N^n(4c)$, we have
\begin{equation}\label{eqn:3.2}
\begin{aligned}
\tilde h^{m^*}_{ij}=&h^{m^*}_{ij}-\frac{n}{n+2}\big(H^{m^*}\delta_{ij}+H^{i^*}\delta_{jm}
+H^{j^*}\delta_{im}\big)\\
=&h^{m^*}_{ij}-c^{m^*}_{ij}
\end{aligned}
\end{equation}
where $c^{m^*}_{ij}=\frac{n}{n+2}\big\{H^{m^*}\delta_{ij}+H^{i^*}\delta_{jm}
+H^{j^*}\delta_{im}\big\}$.

The first covariant derivatives of $\tilde h^{m^*}_{ij}$ are defined by
\begin{equation}\label{eqn:3.3}
\sum_{l=1}^n\tilde h^{m^*}_{ij,l}\theta_l:=d\tilde h_{ij}^{m^*}+\sum_{l=1}^n\tilde h^{m^*}_{il}\theta_{lj}
+\sum_{l=1}^n\tilde h^{m^*}_{jl}\theta_{li}
+\sum_{l=1}^{n}\tilde h^{l^*}_{ij}\theta_{l^*m^*}.
\end{equation}

The second covariant derivatives of $\tilde h^{m}_{ij}$ are defined by
\begin{equation}\label{eqn:3.5}
\sum_{l=1}^n\tilde h^{m^*}_{ij,kl}\theta_l:=d\tilde h_{ij,k}^{m^*}
+\sum_{l=1}^n\tilde h^{m^*}_{lj,k}\theta_{li}
+\sum_{l=1}^n\tilde h^{m^*}_{il,k}\theta_{lj}
+\sum_{l=1}^n\tilde h^{m^*}_{ij,l}\theta_{lk}
+\sum_{l=1}^{n}\tilde h^{\l^*}_{ij,k}\theta_{l^*m^*}.
\end{equation}

On the other hand, we have the following Ricci identities
\begin{equation}\label{eqn:3.5}
\tilde h^{m^*}_{ij,kp}-\tilde h^{m^*}_{ij,pk}
=\sum_{l}\tilde h^{m^*}_{lj}R_{likp}+\sum_l\tilde h^{m^*}_{il}R_{ljkp}
+\sum_{l}\tilde h^{l^*}_{ij}R_{l^*m^*kp}.
\end{equation}

The following proposition links those geometric quantities together:
\begin{lemma}
Let $M^n\hookrightarrow N^n(4c)$ be a Lagrangian submanifold, then the  Lagrangian trace-free second fundamental form $\tilde h$
satisfies
\begin{equation}\label{eqn:3.6}
 | \tilde h |^2= | h |^2-\frac{3n^2}{n+2} | H |^2.
\end{equation}
\begin{equation}\label{eqn:3.7}
\sum_m\tilde h^{m^*}_{ij,m}=\tfrac{n}{n+2}\big(nH^{i^*}_{,j}-{\rm div}JH\,g_{ij}\big).
\end{equation}
\end{lemma}
\begin{proof}
\eqref{eqn:3.6} and \eqref{eqn:3.7} can be immediately obtained from \eqref{eqn:3.2}.
\end{proof}

\begin{definition}[\cite{CU,RU,CMU}]
We define a $(0,2)$-tensor $T$ in local orthonormal basis as follows:
\begin{equation}\label{eqn:3.8}
T_{ij}=\tfrac1n\sum_m\tilde{h}^{m^*}_{ij,m}=\tfrac{1}{n+2}\big(nH^{i^*}_{,j}-\sum_mH^{m^*}_{,m}\,g_{ij}\big)
\end{equation}
\end{definition}
\begin{remark}
$T$ is a trace-free tensor and symmetric. $T=0$ if  and only if $JH$ is a conformal vector field.
\end{remark}

In the following we will derive a Simons' type identity for $\Delta|\tilde h|^2$. First we have
\begin{lemma}\label{lem:3.2}
Let  $M^n\hookrightarrow N^n(4c)$ be a Lagrangian immersion. Then
\begin{equation}\label{eqn:3.9}
\begin{aligned}
\sum_{ijmk}\tilde h^{m^*}_{ij}\tilde h^{m^*}_{ij,kk}=(n+2)\langle\tilde{h},\nabla T\rangle
+ \sum_{i,j,m,k,l}\tilde h^{m^*}_{ij}
\Big(\tilde h^{m^*}_{lk}R_{lijk}
+\tilde h^{m^*}_{il}R_{lkjk}
+\tilde h^{l^*}_{ik}R_{l^*m^*jk}\Big)
\end{aligned}
\end{equation}
\end{lemma}

\begin{proof}
By using the Codazzi equation \eqref{eqn:2.8} and \eqref{eqn:3.2}, the definition of $\tilde h$ under local coordinates is just
\begin{equation}\label{eqn:3.10}
 \tilde h^{m^*}_{ij,k}=\tilde h^{m^*}_{ik,j}+\tfrac{n}{n+2}\big(\delta_{ik}H^{m^*}_{,j}
 +\delta_{km}H^{i^*}_{,j}-\delta_{ij}H^{m^*}_{,k}-\delta_{jm}H^{i^*}_{,k}\big)
\end{equation}

With the help of Ricci identity \eqref{eqn:3.5}, \eqref{eqn:3.8} and \eqref{eqn:3.10}, we have
\begin{equation}\label{eqn:3.11}
\begin{aligned}
 \sum_k \tilde h^{m^*}_{ij,kk}=&\sum_k\tilde h^{m^*}_{ik,jk}+\sum_k\tfrac{n}{n+2}\big(\delta_{ik}H^{m^*}_{,jk}
 +\delta_{km}H^{i}_{,jk}-\delta_{ij}H^{m^*}_{,kk}-\delta_{jm}H^{i^*}_{,kk}\big)\\
 =&\sum_k \tilde h^{m^*}_{ik,kj}+
\sum_{k,l}\tilde h^{m^*}_{lk}R_{lijk}+\sum_{k,l}\tilde h^{m^*}_{il}R_{lkjk}
+\sum_{k,l}\tilde h^{l^*}_{ik}R_{l^*m^*jk}\\
&+\sum_k\tfrac{n}{n+2}\big(\delta_{ik}H^{m^*}_{,jk}
 +\delta_{km}H^{i}_{,jk}-\delta_{ij}H^{m^*}_{,kk}-\delta_{jm}H^{i^*}_{,kk}\big)\\
 =&\sum_k \tilde h^{m^*}_{kk,ij}+
\sum_{k,l}\tilde h^{m^*}_{lk}R_{lijk}+\sum_{k,l}\tilde h^{m^*}_{il}R_{lkjk}
+\sum_{k,l}\tilde h^{l^*}_{ik}R_{l^*m^*jk}\\
&+\sum_k\tfrac{n}{n+2}\big(\delta_{ik}H^{m^*}_{,jk}
 +\delta_{km}H^{i^*}_{,jk}-\delta_{ij}H^{m^*}_{,kk}-\delta_{jm}H^{i^*}_{,kk}\big)+nT_{im,j}.
 \end{aligned}
\end{equation}
Then, by using \eqref{eqn:3.8} and the fact that $\tilde{h}$ is trace free and tri-symmetric, we have
\begin{equation*}
\begin{aligned}
\sum_{i,j,m,k}\tilde h^{m^*}_{ij}\tilde h^{m^*}_{ij,kk}=&\sum_{i,j,m,k}\tilde h^{m^*}_{ij}\Big[\tilde h^{m^*}_{lk}R_{lijk}+\tilde h^{m^*}_{il}R_{lkjk}+\tilde h^{l^*}_{ik}R_{l^*m^*jk}\Big]\\
&+\sum_{m,i,j}\tilde h^{m^*}_{ij}\Big[T_{mj,i}+T_{ij,m}\Big]+n\sum_{m,i,j}\tilde h^{m^*}_{ij}T_{im,j}\\
=&\sum_{i,j,m,k}\tilde h^{m^*}_{ij}\Big[\tilde h^{m^*}_{lk}R_{lijk}+\tilde h^{m^*}_{il}R_{lkjk}+\tilde h^{l^*}_{ik}R_{l^*m^*jk}\Big]\\
&+(n+2)\sum_{m,i,j}\tilde h^{m^*}_{ij}T_{ij,m}
\end{aligned}
\end{equation*}
Thus, we obtain the assertion.
\end{proof}

Next, by using lemma \ref{lem:3.2}
\begin{equation}\label{eqn:3.12}
\begin{aligned}
  \frac12\Delta | \tilde h|^2
=&|\nabla \tilde h|^2+\sum_{ijmk}\tilde h^{m^*}_{ij}\tilde h^{m^*}_{ij,kk}
\\=&|\nabla \tilde h|^2+(n+2)\langle\tilde{h},\nabla T\rangle\\
&+\underbrace{\sum_{i,j,k,m,l}\tilde h^{m^*}_{ij}\tilde h^{m^*}_{lk}R_{lijk}}_{I}
+\underbrace{\sum_{i,j,k,m,l}\tilde h^{m^*}_{ij}\tilde h^{m^*}_{il}R_{lkjk}}_{II}
+\underbrace{\sum_{i,j,k,m,l}\tilde h^{m^*}_{ij}\tilde h^{l^*}_{ik}R_{lmjk}}_{III}.
\end{aligned}
\end{equation}
Note that by the symmetry of $\tilde{h}_{ij}^{k^*}$, $I=III$. Hence we only need to compute $I$ and $II$. Direct computations show that
\begin{equation}\label{eqn:3.13}
\begin{aligned}
I=&c\sum_{i,j,k,m,l}\tilde h^{m^*}_{ij}\tilde h^{m^*}_{kl}(\delta_{lj}\delta_{ik}-\delta_{lk}
  \delta_{ij})+\sum_{i,j,k,m,l,t}\tilde h^{m^*}_{ij}\tilde h^{m^*}_{kl}(\tilde h^{t^*}_{lj}\tilde h^{t^*}_{ik}-\tilde h^{t^*}_{lk}
  \tilde h^{t^*}_{ij})\\
  &+\sum_{i,j,k,m,l,t}\tilde h^{m^*}_{ij}\tilde h^{m^*}_{kl}(\tilde h^{t^*}_{lj}c^{t^*}_{ik}+c^{t^*}_{lj}\tilde h^{t^*}_{ik}-\tilde h^{t^*}_{lk}
  c^{t^*}_{ij}-c^{t^*}_{lk}\tilde h^{t^*}_{ij}+c^{t^*}_{lj}c^{t^*}_{ik}-c^{t^*}_{lk}c^{t^*}_{ij})\\
  =& c| \tilde h |^2+\tfrac{n^2}{(n+2)^2} | \tilde h |^2 | H |^2+\tfrac{2n}{n+2}\sum_{j,k,l,m,t}\tilde h^{m^*}_{jk}
  \tilde h^{m^*}_{kl}\tilde h^{t^*}_{lj}H^{t^*}\\
  &+\sum_{i,j,k,m,l,t}\tilde h^{m^*}_{ij}\tilde h^{m^*}_{kl}(\tilde h^{t^*}_{lj}\tilde h^{t^*}_{ik}-\tilde h^{t^*}_{lk}
  \tilde h^{t^*}_{ij})+\tfrac{2n^2}{(n+2)^2}\sum_{i,j,k,m}\tilde h^{m^*}_{ij}\tilde h^{m^*}_{jk}H^{i^*}H^{k^*},
 \end{aligned}
\end{equation}
where in the second equality we used the following identities derived by direct computations
\begin{eqnarray*}
\sum_{i,j,k,m,l,t}\tilde h^{m^*}_{ij}\tilde h^{m^*}_{kl}\tilde h^{t^*}_{lj}c^{t^*}_{ik}=\sum_{i,j,k,m,l,t}\tilde h^{m^*}_{ij}\tilde h^{m^*}_{kl}c^{t^*}_{lj}\tilde h^{t^*}_{ik}=\tfrac{3n}{n+2}\sum_{j,k,l,m,t}\tilde h^{m^*}_{jk}
  \tilde h^{m^*}_{kl}\tilde h^{t^*}_{lj}H^{t^*},
\end{eqnarray*}
\begin{eqnarray*}
\sum_{i,j,k,m,l,t}\tilde h^{m^*}_{ij}\tilde h^{m^*}_{kl}\tilde h^{t^*}_{lk}
  c^{t^*}_{ij}=\sum_{i,j,k,m,l,t}\tilde h^{m^*}_{ij}\tilde h^{m^*}_{kl}c^{t^*}_{lk}\tilde h^{t^*}_{ij}=\tfrac{2n}{n+2}\sum_{j,k,l,m,t}\tilde h^{m^*}_{jk}
  \tilde h^{m^*}_{kl}\tilde h^{t^*}_{lj}H^{t^*},
\end{eqnarray*}
and since
\begin{equation*}
\begin{aligned}
\sum_tc^{t^*}_{lj}c^{t^*}_{ik}=&\frac{n^2}{(n+2)^2}\sum_t\Big((\mathfrak{S}_{t,l,j}
H^{t^*}\delta_{lj})(\mathfrak{S}_{t,i,k}H^{t^*}\delta_{ik})
\Big)\\
=&\frac{n^2}{(n+2)^2}\Big(\mathfrak{S}_{i,j,k,l}H^{l^*}H^{i^*}\delta_{jk}+2H^{l^*}H^{j^*}\delta_{ik}\\
&+2H^{i^*}H^{k^*}\delta_{jl}+|H|^2\delta_{ik}\delta_{jl}\Big),
\end{aligned}
\end{equation*}
where $\mathfrak{S}$  stands for the cyclic sum, thus
\begin{eqnarray*}
\sum_{i,j,k,m,l,t}\tilde h^{m^*}_{ij}\tilde h^{m^*}_{kl}c^{t^*}_{lj}c^{t^*}_{ik}=\tfrac{n^2}{(n+2)^2} | \tilde h |^2 | H |^2+\tfrac{6n^2}{(n+2)^2}\sum_{i,j,k,m}\tilde h^{m^*}_{ij}\tilde h^{m^*}_{jk}H^{i^*}H^{k^*},
\end{eqnarray*}
\begin{eqnarray*}
\sum_{i,j,k,m,l,t}\tilde h^{m^*}_{ij}\tilde h^{m^*}_{kl}c^{t^*}_{lk}c^{t^*}_{ij}=\tfrac{4n^2}{(n+2)^2}\sum_{i,j,k,m}\tilde h^{m^*}_{ij}\tilde h^{m^*}_{jk}H^{i^*}H^{k^*}.
\end{eqnarray*}
Similarly we have
\begin{equation}\label{eqn:3.14}
\begin{aligned}
II&=\sum_{i,j,k,m,l}\tilde h^{m^*}_{ij}\tilde h^{m^*}_{il}[c(\delta_{lj}\delta_{kk}-\delta_{lk}\delta_{jk})+\sum_t(h_{lj}^{t^*}h_{kk}^{t^*}-h_{lk}^{t^*}h_{kj}^{t^*})
\\&=(n-1)c | \tilde h|^2+\sum_{i,j,k,m,l,t}\tilde h^{m^*}_{ij}\tilde h^{m^*}_{li}\big(n\tilde h^{t^*}_{lj}H^{t^*}+nc^{t^*}_{lj}H^{t^*}\\
&-\tilde h^{t^*}_{lk}
  \tilde h^{t^*}_{kj}-\tilde h^{t^*}_{lk}c^{t^*}_{kj}-c^{t^*}_{lk}\tilde h^{t^*}_{kj}-c^{t^*}_{lk}c^{t^*}_{kj}\big)\\
&=(n-1)c | \tilde h |^2+\tfrac{n^3}{(n+2)^2} | \tilde h |^2 | H |^2+\tfrac{n^2-2n}{n+2}\sum_{ijlmt}\tilde h^{m^*}_{ij}
  \tilde h^{m^*}_{li}\tilde h^{t^*}_{lj}H^{t^*}\\
  &+\tfrac{n^2(n-2)}{(n+2)^2}\sum_{ijml}\tilde h^{m^*}_{ij}\tilde h^{m^*}_{li}H^{j^*}H^{l^*}-\sum_{i,j,k,m,l,t}\tilde h^{m^*}_{ij}\tilde h^{m^*}_{li}\tilde h^{t^*}_{lk}\tilde h^{t^*}_{kj},
 \end{aligned}
\end{equation}
where in the third equality we used the following identities derived by direct computations
\begin{eqnarray*}
n\sum_{i,j,k,m,l,t}\tilde h^{m^*}_{ij}\tilde h^{m^*}_{li}c^{t^*}_{lj}H^{t^*}=\frac{n^2}{n+2} | \tilde h |^2 | H |^2+\frac{2n^2}{n+2}\sum_{ijml}\tilde h^{m^*}_{ij}\tilde h^{m^*}_{li}H^{j^*}H^{l^*},
\end{eqnarray*}
\begin{eqnarray*}
\sum_{i,j,k,m,l,t}\tilde h^{m^*}_{ij}\tilde h^{m^*}_{li}\tilde h^{t^*}_{lk}c^{t^*}_{kj}=\sum_{i,j,k,m,l,t}\tilde h^{m^*}_{ij}\tilde h^{m^*}_{li}\tilde h^{t^*}_{kj}c^{t^*}_{lk}=\frac{2n}{n+2}\sum_{ijlmt}\tilde h^{m^*}_{ij}
  \tilde h^{m^*}_{li}\tilde h^{t^*}_{lj}H^{t^*},
\end{eqnarray*}
and since
\begin{equation*}
\begin{aligned}
\sum_tc^{t^*}_{lk}c^{t^*}_{jk}=&\frac{n^2}{(n+2)^2}\sum_t\Big((\mathfrak{S}_{t,l,k}
H^{t^*}\delta_{lk})(\mathfrak{S}_{t,j,k}H^{t^*}\delta_{jk})
\Big)\\
=&\frac{n^2}{(n+2)^2}\Big(\mathfrak{S}_{j,k,k,l}H^{l^*}H^{j^*}\delta_{kk}+2H^{l^*}H^{k^*}\delta_{jk}\\
&+2H^{j^*}H^{k^*}\delta_{kl}+|H|^2\delta_{jk}\delta_{kl}\Big),
\end{aligned}
\end{equation*}
thus
\begin{eqnarray*}
\sum_{i,j,k,m,l,t}\tilde h^{m^*}_{ij}\tilde h^{m^*}_{li}c^{t^*}_{lk}c^{t^*}_{kj}=\frac{2n^2}{(n+2)^2} | \tilde h |^2 | H |^2+\frac{(n+6)n^2}{(n+2)^2}\sum_{ijml}\tilde h^{m^*}_{ij}\tilde h^{m^*}_{li}H^{j^*}H^{l^*}.
\end{eqnarray*}
\begin{equation}\label{eqn:3.15}
\begin{aligned}
III=I=& c| \tilde h |^2+\tfrac{n^2}{(n+2)^2} | \tilde h |^2 | H |^2+\tfrac{2n}{n+2}\sum_{j,k,l,m,t}\tilde h^{m^*}_{jk}
  \tilde h^{m^*}_{kl}\tilde h^{t^*}_{lj}H^{t^*}\\
  &+\sum_{i,j,k,m,l,t}\tilde h^{m^*}_{ij}\tilde h^{m^*}_{kl}(\tilde h^{t^*}_{lj}\tilde h^{t^*}_{ik}-\tilde h^{t^*}_{lk}
  \tilde h^{t^*}_{ij})\\
  &+\tfrac{2n^2}{(n+2)^2}\sum_{i,j,k,m}\tilde h^{m^*}_{ij}\tilde h^{m^*}_{jk}H^{i^*}H^{k^*}.
 \end{aligned}
\end{equation}

Set $A_{i*}=(\tilde h^{i^*}_{jk})$. Then it follows from \eqref{eqn:3.9} and \eqref{eqn:3.12}-\eqref{eqn:3.15} that
\begin{equation}\label{eqn:3.16}
\begin{aligned}
\frac12\Delta | \tilde h |^2=&(n+2)\langle\tilde{h}, \nabla T\rangle
+|\nabla\tilde h|^2+(n+1) c| \tilde h |^2
  +\tfrac{n^2}{(n+2)} | \tilde h |^2 | H |^2\\
  &+\sum_{i,j}{\rm tr}(A_{i^*}A_{j^*}-A_{j^*}A_{i^*})^2-\sum_{i,j}({\rm tr}A_{i^*}A_{j^*})^2\\
  &+n\sum_{i,j,l,m,t}\tilde h^{m^*}_{ji}\tilde h^{m^*}_{jt}\tilde h^{l^*}_{ti}H^{l^*}
  +\tfrac{n^2}{(n+2)}\sum_{i,j,k,m}\tilde h^{m^*}_{ij}\tilde h^{m^*}_{jk}H^{i^*}H^{k^*}
 \end{aligned}
\end{equation}
Next we estimate terms on the right hand side of (\ref{eqn:3.16}). We will need the following lemma.
\begin{lemma}[\cite{LiLi}]\label{lem:3.3}
Let $B_1,\ldots,B_m$ be symmetric $(n\times n)$-matrices $(m\geq2)$. Denote
$S_{mk}={\rm trace}(B^t_{m}B_{k})$, $S_{m}=S_{mm}=N(B_{m})$,
$S=\sum_{i=1}^mS_{i}$. Then
\begin{equation*}
\sum_{m,k}N(B_{m}B_{k}-B_{k}B_{m})+\sum_{m,k}S^2_{mk}\leq\frac32S^2.
\end{equation*}
\end{lemma}

Now we are prepared to estimate the right hand side of (\ref{eqn:3.16}), mainly the last two terms on the last line of (\ref{eqn:3.16}).

 We re-choose $\{e_i\}_{i=1}^n$ such that $\sum_l\tilde h^{l^*}_{ij}H^{l^*}=\lambda_i\delta_{ij}$, denote by $$S_{H}=\sum_{j,i}(\sum_l\tilde h^{l^*}_{ji}H^{l^*})^2=\sum_{j}\lambda^2_j,\ \ S_{i^*}=\sum_{j,l}(\tilde h^{i^*}_{jl})^2,$$  then $|\tilde{h}|^2=\sum_iS_{i^*}$.
Using lemma \ref{lem:3.3} and the above fact, we have the following estimates for the right hand side of (\ref{eqn:3.16})
\begin{equation}\label{eqn:3.17}
\begin{aligned}
\frac12\Delta | \tilde h |^2\geq&(n+2)\langle\tilde{h},\nabla T\rangle
+|\nabla\tilde h|^2+(n+1)c | \tilde h |^2
  +\tfrac{n^2}{(n+2)} | \tilde h |^2 | H |^2-\tfrac32|\tilde{h}|^4\\
  &+n\sum_i\lambda_iS_{i^*}+\tfrac{n^2}{n+2}\sum_i\lambda^2_{i}\\
  \geq&(n+2)\langle\tilde{h},\nabla T\rangle
+|\nabla\tilde h|^2+(n+1)c | \tilde h |^2
  +\tfrac{n^2}{(n+2)} | \tilde h |^2 | H |^2-\tfrac32|\tilde{h}|^4\\
  &+\tfrac{n}2\sum_i(\lambda_i+S_{i^*})^2-\tfrac{n}2\sum_iS^2_{i^*}\\
  \geq&(n+2)\langle\tilde{h},\nabla T\rangle+|\nabla\tilde h|^2+(n+1)c | \tilde h |^2
  +\tfrac{n^2}{(n+2)} | \tilde h |^2 | H |^2-\tfrac{n+3}2|\tilde{h}|^4\\
  &+\tfrac{n}2\sum_i(|H|\lambda_i+S_{i^*})^2,\\
  \geq&(n+2)\langle\tilde{h},\nabla T\rangle
+|\nabla\tilde h|^2+(n+1)c | \tilde h |^2
  +\tfrac{n^2}{(n+2)} | \tilde h |^2 | H |^2-\tfrac32|\tilde{h}|^4\\
  &+\tfrac{n}2\sum_i(\lambda_i+S_{i^*})^2-\tfrac{n}2\sum_iS^2_{i^*}\\
  \geq&(n+2)\langle\tilde{h},\nabla T\rangle+|\nabla\tilde h|^2+(n+1)c | \tilde h |^2
  +\tfrac{n^2}{(n+2)} | \tilde h |^2 | H |^2-\tfrac{n+3}2|\tilde{h}|^4,\\
 \end{aligned}
\end{equation}
where in the last inequality we have used $|\tilde{h}|^4=(\sum_iS_{i^*})^2\geq\sum_iS_{i^*}^2$.
\section{Proof of Theorem \ref{thm:1.3}}

In this section we will use $C$ to denote constants depending only on $n$, which may vary line by line. We have
\begin{lemma}\label{lem:4.1}
Assume that $M^n$ is a complete Lagrangian submanifold in $\mathbb{C}^n$, and let $\gamma$ be a cut-off function on $M^n$ with $\|\nabla\gamma\|_{L_\infty}=\Gamma,$ then we have
\begin{eqnarray}\label{eqn:4.1}
&&\int_{M}(|\nabla\tilde{h}|^2+|H|^2|\tilde{h}^2|)\gamma^2d\mu\nonumber
\\&& \leq C\int_M\langle\nabla^*T,H\lrcorner\omega\rangle\gamma^2d\mu+C\int_M|\tilde{h}|^4\gamma^2d\mu+C\Gamma^2\int_{\{\gamma>0\}}|h|^2d\mu.
\end{eqnarray}
\end{lemma}
\proof Multiplying (\ref{eqn:3.17}) by $\gamma^2$ we get
\begin{eqnarray}
&&\frac12\int_M\gamma^2\Delta | \tilde h |^2d\mu
\geq(n+2)\int_M\langle\tilde{h}\gamma^2, \nabla T\rangle d\mu\nonumber
\\&&+\int_M|\nabla\tilde h|^2\gamma^2d\mu+\tfrac{n^2}{(n+2)}\int_M | \tilde h |^2 | H |^2\gamma^2d\mu-\tfrac{n+3}2\int_M| \tilde h |^4\gamma^2d\mu.
\end{eqnarray}
Note that
\begin{eqnarray}
\frac12\int_M\gamma^2\Delta | \tilde h |^2d\mu=-\int_M\langle\nabla\tilde{h}, \gamma\nabla\gamma\otimes\tilde{h}\rangle d\mu,
\end{eqnarray}
and by the definition of $T$ and integral by parts
\begin{eqnarray*}
\int_M\langle\tilde{h}\gamma^2, \nabla T\rangle d\mu&=&-\int_M\langle T, \gamma^2\nabla^*\tilde{h}\rangle-2\int_M\langle T\otimes\nabla\gamma,\gamma\tilde{h}\rangle d\mu
\\&=&-n\int_M\langle T,T\rangle\gamma^2d\mu-2\int_M\langle T\otimes\nabla\gamma,\gamma\tilde{h}\rangle  d\mu
\\&=&-n\int_M\langle T,\frac{1}{n+2}(n\nabla H\lrcorner\omega-divJHg)\gamma^2\rangle d\mu-2\int_M\langle T\otimes\nabla\gamma,\gamma\tilde{h}\rangle  d\mu
\\&=&-\frac{n^2}{n+2}\int_M\langle T,\nabla H\lrcorner\omega\rangle\gamma^2\rangle d\mu-2\int_M\langle T\otimes\nabla\gamma,\gamma\tilde{h}\rangle  d\mu
\\&=&\frac{n^2}{n+2}\int_M\langle\nabla^*T,H\lrcorner\omega\rangle\gamma^2d\mu+\frac{2n^2}{n+2}\int_M\langle T,H\lrcorner\omega\otimes\nabla\gamma\rangle\gamma d\mu
\\&-&2\int_M\langle T\otimes\nabla\gamma,\gamma\tilde{h}\rangle  d\mu.
\end{eqnarray*}
Therefore since $|T|\leq C|\nabla\tilde{h}|$, $|\tilde{h}|\leq C|h|$ and $|\nabla\gamma|\leq\Gamma$ we have
\begin{eqnarray*}
&&\int_M|\nabla\tilde h|^2\gamma^2+\tfrac{n^2}{(n+2)} | \tilde h |^2 | H |^2\gamma^2d\mu
\\&\leq&-\int_M\langle\nabla\tilde{h}, \gamma\nabla\gamma\otimes\tilde{h}\rangle d\mu-n^2\int_M\langle\nabla^*T,H\lrcorner\omega\rangle\gamma^2d\mu-2n^2\int_M\langle T,H\lrcorner\omega\otimes\nabla\gamma\rangle\gamma d\mu
\\&+&2(n+2)\int_M\langle T\otimes\nabla\gamma,\gamma\tilde{h}\rangle  d\mu+\tfrac{n+3}2\int_M| \tilde h |^4\gamma^2d\mu
\\&\leq&-n^2\int_M\langle\nabla^*T,H\lrcorner\omega\rangle\gamma^2d\mu+\tfrac{n+3}2\int_M| \tilde h |^4\gamma^2d\mu+\frac12\int_M|\nabla\tilde h|^2\gamma^2d\mu+C\Gamma^2\int_{\{\gamma>0\}}|h|^2d\mu.
\end{eqnarray*}
\endproof

We will need the following Michael-Simon inequality:
\begin{theorem}[\cite{MS}\cite{HS}]
Assume that $M^n$ is a compact submanifold of $R^{n+p}$ with or without boundary. Assume that $v \in C^1(M^n)$ is a nonnegative function such that $v=0$ on $\partial M^n$, if $\partial M^n$ is not an empty set. Then
\begin{eqnarray*}
(\int_Mv^{\frac{n}{n-1}}d\mu)^{\frac{n-1}{n}}\leq C\int_M|\nabla v|+v|H|d\mu,
\end{eqnarray*}
where $H$ is the mean curvature vector field of $M^n$.
\end{theorem}
Now assume that $n\geq3$. In the above Michael-Simon inequality we let $v=f^{\frac{2(n-1)}{n-2}}$, then by H\"older inequality we easily get
\begin{eqnarray}\label{eqn:3.20}
(\int_Mf^{\frac{2n}{n-2}}d\mu)^{\frac{n-2}{n}}\leq C\int_M|\nabla f|^2+f^2|H|^2d\mu.
\end{eqnarray}
Therefore by letting $f=|\tilde{h}|\gamma$ in (\ref{eqn:3.20}) we obtain
\begin{eqnarray}\label{eqn:3.21}
\int_M|\tilde{h}|^4\gamma^2d\mu\nonumber
&\leq&(\int_M|\tilde{h}|^nd\mu)^{\frac2n}(\int_M(|\tilde{h}|\gamma)^{\frac{2n}{n-2}}d\mu)^{\frac{n-2}{n}}\nonumber
\\&\leq&C(\int_M|\tilde{h}|^nd\mu)^{\frac2n}(\int_M|\nabla(|\tilde{h}|\gamma)|^2+|H|^2|\tilde{h}|^2\gamma^2d\mu)\nonumber
\\&\leq& C(\int_M|\tilde{h}|^nd\mu)^{\frac2n}(\int_M|\nabla\tilde{h}|^2\gamma^2d\mu
+\Gamma^2\int_{\{\gamma>0\}}|\tilde{h}|^2d\mu+\int_M|H|^2|\tilde{h}|^2\gamma^2d\mu)\nonumber
\\&\leq&C(\int_M|\tilde{h}|^nd\mu)^{\frac2n}(|\nabla\tilde{h}|^2\gamma^2+|H|^2|\tilde{h}|^2\gamma^2d\mu)\nonumber
\\&+&C\Gamma^2(\int_M|\tilde{h}|^nd\mu)^{\frac2n}\int_{\{\gamma>0\}}|h|^2d\mu.
\end{eqnarray}
From (\ref{eqn:3.17}) and (\ref{eqn:3.21}) we have
\begin{eqnarray*}
&&\int_{M}(|\nabla\tilde{h}|^2+|H|^2|\tilde{h}^2|)\gamma^2d\mu\nonumber
\\&& \leq C\int_M\langle\nabla^*T,H\lrcorner\omega\rangle\gamma^2d\mu+C(\int_M|\tilde{h}|^nd\mu)^{\frac2n}(|\nabla\tilde{h}|^2\gamma^2
+|H|^2|\tilde{h}|^2\gamma^2d\mu)+C\Gamma^2\int_{\{\gamma>0\}}|h|^2d\mu,
\end{eqnarray*}
which implies that if there exists $\epsilon_0$ sufficiently small such that $$\int_M|\tilde{h}|^nd\mu\leq\epsilon_0,$$ we have
\begin{eqnarray}
\int_{M}(|\nabla\tilde{h}|^2+|H|^2|\tilde{h}^2|)\gamma^2d\mu\leq C\int_M\langle\nabla^*T,H\lrcorner\omega\rangle\gamma^2d\mu+C\Gamma^2\int_{\{\gamma>0\}}|h|^2d\mu.
\end{eqnarray}

\textbf{Case 1:} If $\nabla^*T=0$, note that for any $R>0$ we can choose $\gamma\in C_c^1(M_R(p_0))$ such that $\gamma=1$ on $M_\frac R2(p_0)$ where $M_r(p_0)$ denotes geodesic ball of radius $r$ with center $p_0\in M^n$, and $\Gamma\leq\frac{C}{R}$,
therefore by letting $R\to+\infty$ we get $$\int_{M}|\nabla\tilde{h}|^2+|H|^2|\tilde{h}^2|d\mu=0,$$
which implies that $\tilde{h}=0$ and $M^n$ is either a Lagrangian subspace or a Whitney sphere by \cite{CU,RU} or $H=0, \nabla h=0$ and $M^n$ is a Lagrangian subspace by \cite{La}.

\textbf{Case 2:} If $\nabla^*\nabla^*T=0$ and $M^n$ is a Lagrangian sphere, then by Dazord \cite{Da}, there exists a smooth function $f$
on $M^n$ such that $H\lrcorner\omega=df$, and let $\gamma\equiv1$ on $M^n$, then we have
\begin{eqnarray*}
\int_{M}|\nabla\tilde{h}|^2+|H|^2|\tilde{h}^2|d\mu&\leq& C\int_M\langle\nabla^*T,df \rangle d\mu
\\&=&-C\int_Mf\nabla^*\nabla^*Td\mu
\\&=&0,
\end{eqnarray*}
which implies that $\tilde{h}=0$ and $M^n$ is either a Lagrangian subspace or a Whitney sphere by \cite{CU,RU} or $H=0, \nabla h=0$ and $M^n$ is a Lagrangian subspace by \cite{La}.

This completes the proof of Theorem \ref{thm:1.3}.

\section{Proof of Theorem \ref{thm:1.4}}
The proof of Theorem \ref{thm:1.4} is quite similar with the proofs of Theorem 1.2 in \cite{Zh}, Theorem 1.5 in \cite{LY}
and Theorem \ref{thm:1.3} in the present paper. Therefore we will only give a outline of the proof and omit some details. In this section we will use $C$ to denote constants depending only on $n$ which may vary line by line.

Letting $c=1$ in (\ref{eqn:3.17}), we have
\begin{eqnarray}\label{eqn:5.1}
\frac12\Delta | \tilde h |^2\geq(n+2)\langle\tilde{h},\nabla T\rangle+|\nabla\tilde h|^2+(n+1) | \tilde h |^2
  +\tfrac{n^2}{(n+2)} | \tilde h |^2 | H |^2-\tfrac{n+3}2|\tilde{h}|^4.
\end{eqnarray}

Then similarly with Lemma \ref{lem:4.1} we can obtain the following Lemma.
\begin{lemma}
Assume that $M^n$ is a complete Lagrangian submanifold in $\mathbb{CP}^n$ and $\gamma$ is a cut off function on $M^n$ with $\|\nabla\gamma\|_{L_\infty}=\Gamma,$ then we have
\begin{eqnarray}\label{eqn:5.2}
&&\int_{M}(|\nabla\tilde{h}|^2+|\tilde{h}|^2|H|^2+|\tilde{h}|^2)\gamma^2d\mu\nonumber
\\&& \leq C\int_M\langle\nabla^*T,H\lrcorner\omega\rangle \gamma^2d\mu+C\int_M|\tilde{h}|^4\gamma^2d\mu+C\Gamma^2\int_{\{\gamma>0\}}|h|^2d\mu.
\end{eqnarray}
\end{lemma}
The same as the previous section, to absorb the "bad term" $\int_M|\tilde{h}|^4\gamma^2d\mu$ on the right hand side of (\ref{eqn:5.2}), we will use the Michael-Simom inequality. In order to do this, we need isometrically immersed $\mathbb{CP}^n$ into some Euclidean space $\mathbb{R}^{n+p}$, which is possible by Nash's celebrated embedding theorem. Assume that $\mathbb{CP}^n$ has mean curvature $H_0$ as a submanifold in $\mathbb{R}^{n+p}$, and $M^n$ has mean curvature $\bar{H}$ as a submanifold in $\mathbb{R}^{n+p}$. Then it is easy to see that $|\bar{H}|^2\leq|H_0|^2+|H|^2$.

If $n=2$, from the original Michael-Simon inequality we see that if $M\hookrightarrow \mathbb{R}^{n+p}$ is compact with or without boundary then
\begin{eqnarray}
\int_Mf^2d\mu\leq C\big(\int_M|\nabla f|d\mu+\int_Mf|H|d\mu\big)^2,
\end{eqnarray}
for any nonnegative function $f\in C^1(M)$ with $f|_{\partial M}=0$. Let $f=|\tilde{h}|^2\gamma$ in the above inequality we obtain
\begin{eqnarray}\label{eqn:5.4}
\int_M|\tilde{h}|^4\gamma^2d\mu&\leq& C\big(\int_M|\tilde{h}||\nabla\tilde{h}|\gamma+|\nabla\gamma||\tilde{h}|^2 d\mu+\int_M|\tilde{h}|^2|\bar{H}|\gamma d\mu\big)^2\nonumber
\\&\leq&C\int_M|\tilde{h}|^2d\mu\big(\int_M|\nabla\tilde{h}|^2\gamma^2d\mu+\int_M|\tilde{h}|^2|\bar{H}|^2\gamma^2d\mu\big)
+C\Gamma^2(\int_{\{\gamma>0\}}|\tilde{h}|^2d\mu)^2\nonumber
\\&\leq&C\int_M|\tilde{h}|^2d\mu\big(\int_M|\nabla\tilde{h}|^2\gamma^2d\mu+\int_M|\tilde{h}|^2|H|^2\gamma^2d\mu\big)\nonumber
\\&+&C\max_{\mathbb{CP}^n}|H_0|^2(\int_M|\tilde{h}|^2\gamma^2d\mu)^2+C\Gamma^2(\int_{\{\gamma>0\}}|\tilde{h}|^2d\mu)^2 \nonumber
\\&=&C\int_M|\tilde{h}|^2d\mu\big(\int_M|\nabla\tilde{h}|^2\gamma^2d\mu+\int_M|\tilde{h}|^2|H|^2\gamma^2d\mu\big)\nonumber
\\&+&C(\int_M|\tilde{h}|^2\gamma^2d\mu)^2+C\Gamma^2(\int_{\{\gamma>0\}}|\tilde{h}|^2d\mu)^2.
\end{eqnarray}
From (\ref{eqn:5.2}) and (\ref{eqn:5.4}) we see that
\begin{eqnarray*}
&&\int_{M}(|\nabla\tilde{h}|^2+|\tilde{h}|^2|H|^2|+|\tilde{h}|^2)\gamma^2d\mu
\\&\leq& C\int_M\langle\nabla^*T,H\lrcorner\omega\rangle\gamma^2 d\mu
+C\int_M|\tilde{h}|^2d\mu\big(\int_M|\nabla\tilde{h}|^2\gamma^2d\mu+\int_M|\tilde{h}|^2|H|^2\gamma^2d\mu\big)\nonumber
\\&+&C(\int_M|\tilde{h}|^2\gamma^2d\mu)^2+C\Gamma^2(\int_{\{\gamma>0\}}|\tilde{h}|^2d\mu)^2+C\Gamma^2\int_{\{\gamma>0\}}|h|^2d\mu.
\end{eqnarray*}
Therefore if $M^n$ satisfies assumptions of Theorem \ref{thm:1.4} we can similarly with the previous section obtain that $\tilde{h}=0$ and hence $M^n$ is the real projective space $\mathbb{RP}^n$ or a Whitney sphere, by \cite{Ch}.

If $n\geq3$, similarly with the proof of Theorem \ref{thm:1.3}, we can obtain from the Michael-Simon inequality that
\begin{eqnarray}
\big(\int_Mf^\frac{2n}{n-2}d\mu\big)^\frac{n-2}{n}\leq C\big(\int_M|\nabla f|^2d\mu+\int_Mf^2|\bar{H}|^2d\mu\big),
\end{eqnarray}
for any nonnegative function $f\in C^1(M^n)$ with $f|_{\partial M}=0$. Therefore by letting $f=|\tilde{h}|\gamma$ in the above inequality we obtain
\begin{eqnarray}\label{eqn:5.6}
\int_M|\tilde{h}|^4\gamma^2d\mu\nonumber
&\leq&(\int_M|\tilde{h}|^nd\mu)^{\frac2n}(\int_M(|\tilde{h}|\gamma)^{\frac{2n}{n-2}}d\mu)^{\frac{n-2}{n}}\nonumber
\\&\leq&C(\int_M|\tilde{h}|^nd\mu)^{\frac2n}(\int_M|\nabla(|\tilde{h}|\gamma)|^2+|\bar{H}|^2|\tilde{h}|^2\gamma^2d\mu)\nonumber
\\&\leq& C(\int_M|\tilde{h}|^nd\mu)^{\frac2n}(\int_M|\nabla\tilde{h}|^2\gamma^2d\mu
+\Gamma^2\int_{\{\gamma>0\}}|\tilde{h}|^2d\mu+\int_M|\bar{H}|^2|\tilde{h}|^2\gamma^2d\mu)
\nonumber
\\&\leq&C(\int_M|\tilde{h}|^nd\mu)^{\frac2n}(\int_M|\nabla\tilde{h}|^2\gamma^2+|H|^2|\tilde{h}|^2\gamma^2d\mu)\nonumber
\\&+&C\max_{\mathbb{CP}^n}|H_0|^2(\int_M|\tilde{h}|^nd\mu)^{\frac2n}\int_M|\tilde{h}|^2\gamma^2d\mu\nonumber
+C\Gamma^2(\int_M|\tilde{h}|^nd\mu)^{\frac2n}\int_{\{\gamma>0\}}|\tilde{h}|^2d\mu\nonumber
\\&=&C(\int_M|\tilde{h}|^nd\mu)^{\frac2n}(\int_M|\nabla\tilde{h}|^2\gamma^2+|H|^2|\tilde{h}|^2\gamma^2d\mu)
+C(\int_M|\tilde{h}|^nd\mu)^{\frac2n}\int_M|\tilde{h}|^2\gamma^2d\mu\nonumber
\\&+&C\Gamma^2(\int_M|\tilde{h}|^nd\mu)^{\frac2n}\int_{\{\gamma>0\}}|\tilde{h}|^2d\mu.
\end{eqnarray}
From (\ref{eqn:5.2}) and (\ref{eqn:5.6}) we obtain
\begin{eqnarray*}
&&\int_{M}(|\nabla\tilde{h}|^2+|\tilde{h}|^2|H|^2+|\tilde{h}|^2)\gamma^2d\mu\nonumber
\\&& \leq C\int_M\langle\nabla^*T,H\lrcorner\omega\rangle\gamma^2 d\mu+C(\int_M|\tilde{h}|^nd\mu)^{\frac2n}\int_M(|\nabla\tilde{h}|^2+|\tilde{h}|^2|H|^2)\gamma^2d\mu\nonumber
\\&+&C(\int_M|\tilde{h}|^nd\mu)^{\frac2n}\int_M|\tilde{h}|^2\gamma^2d\mu
+C\Gamma^2(\int_M|\tilde{h}|^nd\mu)^{\frac2n}\int_{\{\gamma>0\}}|h|^2d\mu+C\Gamma^2\int_{\{\gamma>0\}}|h|^2d\mu.
\end{eqnarray*}
Therefore if $M^n$ satisfies assumptions of Theorem \ref{thm:1.4} we can similarly with the previous section obtain that $\tilde{h}=0$ and hence $M^n$ is the real projective space $\mathbb{RP}^n$ or a Whitney sphere, by \cite{Ch}.

This completes the proof of Theorem \ref{thm:1.4}.

\end{document}